\documentstyle[11pt,leqno]{article}
\newtheorem{theorem}{{\sc Theorem}}
\newcommand{\bt}{\begin{theorem}}
\newcommand{\et}{\end{theorem}}
\newcommand{\newsection}[1]{\setcounter{equation}{0} \setcounter{theorem}{0}
\section{#1}}

\newcommand{\NI}{\noindent}
\newcommand{\bea}{\begin{eqnarray}}
\newcommand{\eea}{\end{eqnarray}}

\def \spec#1 {\mathop{#1}}

\def \b #1 {\bf #1}

\newcommand {\CC}{\centerline}

\newcommand{\ity}{\infty}
\newcommand{\raro}{\rightarrow}

\newcommand{\vsp}{\vskip 1em}
\newcommand{\vspp}{\vskip 2em}

\newcommand{\be}{\begin{equation}}
\newcommand{\ee}{\end{equation}}
\newcommand{\ben}{\begin{eqnarray*}}
\newcommand{\een}{\end{eqnarray*}}

\begin{document}

\sloppy
\CC {\bf Parameter Estimation for Generalized Mixed Fractional }
\CC {\bf Stochastic Heat Equation}
\vsp
\CC {B.L.S. PRAKASA RAO}
\CC {CR Rao Advanced Institute of Mathematics, Statistics}
\CC{and Computer Science, Hyderabad, India}
\vspp
\vspp
\NI{\bf Abstract:} We study the properties of a stochastic heat equation with a generalized mixed fractional Brownian  noise. We obtain the covariance structure, stationarity and obtain bounds for the asymptotic behavior of the solution. We suggest  estimators for the unknown parameters based on discrete time observations and study their asymptotic properties.
\vsp
\NI{\bf Keywords}: Stochastic  partial differential equation ; Generalized mixed fractional Brownian motion; Parameter estimation.
\vsp
\NI {\bf AMS Subject Classification (2020)}: 60G22.

\newsection{Introduction} Generalized mixed fractional Brownian motion is a finite mixture of independent fractional Brownian motions. It is known that this process is a centered Gaussian process which is self-similar in a suitable sense and not Markov. Properties of these processes are studied in Chapter 3 of Mishura and Zili (2018). Long range dependence property of time-changed generlized mixed fractional Brownian motion is investigated in Prakasa Rao (2025). Ergodic properties of the solution of a fractional stochastic equation driven by a mixed fractional Brownian motion are discussed in Avetisian and Ralchenko (2020). Parameter estimation for a mixed fractional stochastic heat equation has been investigated in Avetisian and Ralchenko (2023). Our aim in this paper is to study problems of parameter estimation in a stochastic heat equation driven by a generalized mixed fractional Brownian motion  based on discrete time observations. Our techniques are analogous to those in Avetisian and Ralchenko (2023).
\vsp
We study parameter estimation for a stochastic heat equation of the form
\be
(\frac{\partial u}{\partial t}-\frac{1}{2}\frac{\partial^2u}{\partial x^2})(t,x)= \sigma_1 \dot W^{H_1}_x+\sigma_2 \dot W^{H_2}_x, t >0, x \in R
\ee
under the condition 
\be
u(0,x)=0, x \in R.
\ee
The process on the right hand side of the equation (1.1) is a generalized mixed fractional noise. It consists of two independent fractional Brownian motions $W^{H_1}$ and  $W^{H_2}$ with Hurst indices $H_1,H_2$ and the parameters  $\sigma_1$ and $\sigma_2$ are positive constants. We investigate the problem of estimating the parameters $H_1,H_2,\sigma_1$ and $\sigma_2$ based on discrete observations of the solution $u(t,x)$ of the equation (1.1).
\vsp
\newsection{Preliminaries} 
Suppose that $W^{H_i}=\{W^{H_i}_x, x \in R\}, i=1,2$ are two independent two-sided fractional Brownian motions with Hurst indices $H_i, i=1,2$ respectively. Let $G$ be the Green's function of the heat equation given by
\ben 
G(t,x)&=& \frac{1}{\sqrt{2\pi t}}\exp\{-\frac{x^2}{t}\},\; \mbox{if}\;t > 0,\\ 
&=& \delta_0(x)\; \mbox{if}\; t=0.\\
\een
The random field $\{u(t,x), t \geq 0, x \in R\}$ defined by
\be
u(t,x)=\sigma_1\int_0^t\int_R G(t-s,x-y)dW^{H_1}_y ds+ \sigma_2\int_0^t\int_R G(t-s,x-y)dW^{H_2}_y ds
\ee
is called a solution of the stochastic partial differential equation (SPDE) defined by (1.1) and (1.2). 
\vsp
As pointed out by Avetisian and Ralchenko (2020), the stochastic integrals in (2.1) exist as pathwise Riemann-Stieltjes integrals since the Green's function is Lipshitz continuous and the sample paths of the fractional Brownian motion $W^{H_i}$ are Holder continuous  up to order 
$H_i, i=1,2.$
We will now derive some properties of the solution $u(t,x).$
\vsp
\NI{\bf Theorem 2.1} Let $u(t,x), t\in [0,T], x \in R\}$ be a solution to the equation (1.1) and (1.2) as defined by (2.1). Then the following properties hold.\\
(i)For $0\leq t,s\leq T$ and for $x,z \in R,$
\bea
\;\;\;\;\\\nonumber
Cov(u(t,z), u(s,x+z))& =& Cov(u(t,0),u(s,x))\\\nonumber
&=& \frac{1}{\sqrt{2\pi}}\int_0^t\int_0^s (q+r)^{-\frac{3}{2}}\int_R(\sigma_1^2H_1|y|^{2H_1-1}+\sigma_2^2H_2|y|^{2H_2-1})\\\nonumber 
&&\;\;\;\;\times (sign\; y)(y-x)\exp\{-\frac{(y-x)^2}{2(q+r)}\} dy\;dq\;dr.\\\nonumber
\eea
(ii) For any fixed $t_1,\dots,t_n \in [0,T]$, the multivariate  process $\{(u(t_1,x),\dots,u(t_n,x)), x \in R\}$ is a centered stationary Gaussian process.\\
(iii) The variance of $u(t,x)$ is given by
\be
Var(u(t,x))=E[u(t,x)]^2= \sigma_1^2 v_t(H_1)+\sigma_2^2v_t(H_2), t >0, x \in R
\ee
where
\be
v_t(H)= c_Ht^{H+1} \;\;\mbox{and}\;\;c_H= \frac{2^{H+1}(2^H-1)\Gamma(H+\frac{1}{2})}{\sqrt{\pi}(H+1)}.
\ee
(iv)For $t,s \in [0,T]$ and $x>0,$ the covariance function admits the following upper bound:
\be 
|Cov(u(t,0),u(s,x))|\leq C_{H_1,H_2} ts (\sigma_1^2x^{2H_1-2}+\sigma_2^2x^{2H_2-2})
\ee
where $C_{H_1,H_2}$ is a positive constant depending on $H_1$ and $H_2.$\\
(v) For $t,s \in [0,T]$ and $x\in R,$
\bea
Cov(u(t,x),u(s,x))&=& \frac{\sigma_1^22^{H_1}\Gamma(H_1+\frac{1}{2})((t+s)^{H_1+1}-t^{H_1+1}-s^{H_1+1})}{\sqrt{\pi}(H_1+1)}\\\nonumber
&&\;\;\;+ \frac{\sigma_2^22^{H_2}\Gamma(H_2+\frac{1}{2})((t+s)^{H_2+1}-t^{H_2+1}-s^{H_2+1})}{\sqrt{\pi}(H_2+1)}.\\\nonumber
\eea
(vi) For fixed $t>0,$ the process $\{u(t,x), x \in R\}$ is ergodic.
\vsp
\NI{\bf Proof:} Since the processes $W^{H_1}$ and $W^{H_2}$ are independent fractional Brownian motions, it is easy to see that
\be
Cov(u(t,x),u(s,z))= \sigma_1^2 Cov(u_1(t,x),u_1(s,z))+\sigma_2^2 Cov(u_2(t,x),u_2(s,z))
\ee
where 
$$u_i(t,x)=\int_0^t\int_RG(t-s,x-y) dW_y^{H_i} ds, i=1,2.$$
The properties (i)-(v) follow from the results in Avetisian and Ralchenko (2020,2023) for fractional Brownian motion as noise. For any fixed $t\in [0,T],$ the process $\{u(t,x), x \in R\}$ is a stationary Gaussian process . From Proposition 4 in Avetisian and Ralchenko (2020), it follows that the covariance function $R(t,x)=Cov(u(t,0),u(t,x))$ of the process will satisfy the inequality
$$|R(t,x)|\leq C_{H_1}\sigma_1^2t^2x^{2H_1-2}+  C_{H_2}\sigma_1^2t^2x^{2H_2-2}, x>0.$$  
Since $0< H_1,H_2<1,$, it follows that $R(t,x) \raro 0$ as $x \raro \infty.$ This in turn implies that the process $\{u(t,x), x \in R\}$ is an ergodic process for any fixed $t >0.$
\vsp
Let $\delta > 0$ and define
\be
V_N(t)= \frac{1}{N} \sum_{i=1}^N [u(t,k\delta)]^2, t >0, N \geq 1,
\ee
\be
\mu(t)= \sigma_1^2v_t(H_1)+\sigma_2^2v_t(H_2).
\ee
Let
$$
\rho_{t,s}^{H_1,H_2}(k)=Cov(u(t,k\delta),u(s,0)) \;\;\mbox{and}\;\; r_{t,s}(H_1,H_2)=2\sum_{k=-\ity}^{\ity}[\rho_{t,s}^{H_1,H_2}(k)]^2.
$$
\vsp
\NI{\bf Theorem 2.2} (i) For any $t>0,$
\be
V_N(t) \raro \mu(t) a.s.\;\; \mbox{as} \;\; N\raro \infty.
\ee
(ii)Suppose further that $H_1,H_2 \in (0,\frac{3}{4}). $ Then, for any distinct positive $t_1,\dots,t_n,$ the random vector 
$$\sqrt{N}(V_N(t_1)-\mu(t_1),\dots, V_N(t_n)-\mu(t_n))$$
converges in law to a multivariate normal distribution with mean vector 0 and covariance matrix $R$
where 
$$R=((r_{t_i,t_j}(H_1,H_2))_{n\times n}.$$
\vsp
\NI{\bf Proof:} Since the process $\{u(t,x), x \in R\}$ is an ergodic process for any $t >0,$, it follows that 

$$V_N(t)=\frac{1}{N}\sum_{k=1}^N [u(t,k\delta)]^2 \raro E([u(t,0)]^2)\;\; a.s \mbox{as}\;\; N \raro \ity.$$

Note that
$$|\rho_{t_i,t_j}^{H_1,H_2}(k)|\leq C(\sigma_1^2 (k\delta)^{2H_1-2}+\sigma_2^2 (k\delta)^{2H_2-2})$$
for some constant $C>0$ depending on $H_1$ and $H_2$ and hence
$$ (\rho_{t_i,t_j}^{H_1,H_2}(k))^2\leq C (\sigma_1^4 (k\delta)^{4H_1-4}+\sigma_2^4 (k\delta)^{4H_2-4})$$
which in turn implies that 
$$\sum_{k=-\ity}^{\ity}(\rho_{t_i,t_j}^{H_1,H_2}(k))^2<\ity$$
since $H_1,H_2\in (0,\frac{3}{4}).$ Following the arguments in the proof of Theorem 1 in Avetisian and Ralchenko (2023), it follows that 
$$\sqrt{N}(V_N(t_1)-\mu(t_1),\dots, V_N(t_n)-\mu(t_n))$$
converges in law to a multivariate normal distribution with mean vector 0 and covariance matrix $R$
where 
$$R=((r_{t_i,t_j}(H_1,H_2))_{n\times n}$$
by the Cramer-Wold technique and the multivariate Breuer-Major theorem (cf. Arcones (1994)).
\vsp
\newsection{Estimation of $H_1$ given $H_2.$  } 
Let $\delta >0.$ We now consider the problem of estimation of the parameter $H_2$ given $H_1, \sigma_1^2$ and $\sigma_2^2$ and the process $\{u(t,x), t \geq 0, x\in r\}$ is observed at the points, $x_k=k \delta, k=1,\dots,N$ for fixed $t_1,\dots,t_n.$ Following the method of moments for estimation of the unknown parameters which consists in equating the sample moments to the population moments and observing that 
$$ V_N(t) \raro \mu(t)= \sigma_1^2v_t(H_1)+\sigma_2^2v_t(H_2)\;\; a.s \;\;\mbox{as} \; \; N\raro \ity,$$
we obtain the equations
\bea
V_n(t_i)&=& \sigma_1^2 v_{t_i}(H_1)+\sigma_2^2v_{t_i}(H_2), i=1,2,3\\\nonumber
&=& \sigma_1^2 C_{H_1}t_i^{H_1+1}+\sigma_2^2 C_{H_2}t_i^{H_2+1}, i=1,2,3.\\\nonumber
\eea
As a consequence, it follows that
\be
t_2^{-(H_2+1)}V_N(t_2)-t_1^{-(H_2+1)}V_N(t_1)=\sigma_1^2C_{H_1}(t_2^{H_1-H_2}-t_1^{H_1-H_2}) 
\ee
and
\be
t_3^{-(H_2+1)}V_N(t_3)-t_1^{-(H_2+1)}V_N(t_1)=\sigma_1^2C_{H_1}(t_3^{H_1-H_2}-t_1^{H_1-H_2}).
\ee
Taking ratios of the terms on either side of the above equations, we obtain that
$$\frac{t_2^{-(H_2+1)}V_N(t_2)-t_1^{-(H_2+1)}V_N(t_1)}{t_3^{-(H_2+1)}V_N(t_3)-t_1^{-(H_2+1)}V_N(t_1)}= \frac{(t_2^{H_1-H_2}-t_1^{H_1-H_2)}}{(t_3^{H_1-H_2}-t_1^{H_1-H_2})}.$$
Observe that
$$\lim_{H_1\raro H_2}\frac{(t_2^{H_1-H_2}-t_1^{H_1-H_2})}{(t_3^{H_1-H_2}-t_1^{H_1-H_2})}= \frac{\log t_2-\log t_1}{\log t_3 -\log t_1}$$
by L'Hopital rule.
Define the function
\bea
f(H) &= &\frac{(t_2^{H-H_2}-t_1^{H-H_2})}{(t_3^{H-H_2}-t_1^{H-H_2})}\;\;\mbox{if}\;\; H\neq H_2\\\nonumber
&= & \frac{\log t_2-\log t_1}{\log t_3 -\log t_1} \;\;\mbox{if}\;\; H=H_2.\\\nonumber
\eea
For any fixed $H_2,$ and for $0<t_1<t_2<t_3,$ it can be shown that the function $f: R\raro (0,\ity)$ is strictly increasing function in $H$ following arguments in Avetisian and Ralchenko (2023) and hence has an inverse  $f^{-1}.$  We define the estimator $\hat H_{1N}$ of the parameter $H_1$ by the equation 
\be
\hat H_{1N} = f^{-1}(\frac{t_2^{-(H_2+1)}V_N(t_2)-t_1^{-(H_2+1)}V_N(t_1)}{t_3^{-(H_2+1)}V_N(t_3)-t_1^{-(H_2+1)}V_N(t_1)})
\ee
which will be well-defined for large $N.$ Following the method of proof of Theorem 1 in Avetisian and Ralchenko (2023), we obtain the following result.
\vsp
\NI{\bf Theorem 3.1} Suppose $H_1\in (0,1)$ and $H_1\neq H_2.$ Then the estimator $\hat H_{1N}$ is a strongly consistent estimator of $H_1$ as $N\raro \ity.$ Furthermore 
$$\sqrt{N} (\hat H_{1N}-H_1) \raro N(0,\zeta^2) \;\;\mbox{in distribution as}\;\;N\raro \ity$$
where $\zeta^2$ depends on $t_1,t_2,t_3$ and $H_2.$
\vsp
\NI{\bf Proof:} Note that
$$\frac{t_2^{-(H_2+1)}V_N(t_2)-t_1^{-(H_2+1)}V_N(t_1)}{t_3^{-(H_2+1)}V_N(t_3)-t_1^{-(H_2+1)}V_N(t_1)} \raro \;\;f(H_1)\;\; a.s \;\;\mbox{as} \;\;N\raro \ity.$$
From the continuity of the inverse function $f^{-1},$ it follows that $\hat H_{1N}$ converges a.s. to $H_1$ as $N \raro \ity.$ Taking expectations on both sides of the equations (3.2) and (3.3) and then taking the ratios, we obtain that
$$\frac{t_2^{-(H_2+1)}\mu(t_2)-t_1^{-(H_2+1)}\mu(t_1)}{t_3^{-(H_2+1)}\mu(t_3)-t_1^{-(H_2+1)}\mu(t_1)}=\frac{(t_2^{H_1-H_2}-t_1^{H_1-H_2})}{(t_3^{H_1-H_2}-t_1^{H_1-H_2})}=f(H_1).$$
Hence
$$H_1=f^{-1}(\frac{t_2^{-(H_2+1)}\mu(t_2)-t_1^{-(H_2+1)}\mu(t_1)}{t_3^{-(H_2+1)}\mu(t_3)-t_1^{-(H_2+1)}\mu(t_1)}).$$
Therefore
$$\sqrt {N} (\hat H_{1N}-H_1)= \sqrt{N}(g(V_n(t_1),V_N(t_2),V_n(t_3))-g(\mu(t_1),\mu(t_2),\mu(t_3)))$$
where
$$g(x_1,x_2,x_3)= f^{-1}(\frac{t_2^{-(H_2+1)}x_2-t_1^{-(H_2+1)}x_1}{t_3^{-(H_2+1)}x_3-t_1^{-(H_2+1)}x_1}).$$
Applying the delta method and observing that $(V_N(t_1),V_N(t_2),V_N(t_3))$ is asymptotically normal after suitable scaling, it can be shown that 
$$\sqrt{N} (\hat H_{1N}-H_1) \raro N(0,\zeta^2) \;\;\mbox{in distribution as}\;\;N\raro \ity$$
for some $\zeta^2$ depending on $T_i, i=1,2,3$ and $H_2.$. We skip the details.
\vsp
\newsection{Estimation of $\sigma_1^2$ and $\sigma_2^2$ when $H_1$ and $H_2$ are known and $H_1 \neq H_2$}
Suppose that the Hurst indices $H_1$ and $H_2$ are known. We now study the problem of estimation of the parameters $\sigma_1^2$ and $\sigma_2^2$ based on the discrete set of observations $u(t_i,k\delta), i=1,2,, k=1,\dots,N$ with $t_1<t_2$ and a fixed $\delta >0.$ Using the method of moments again, we obtain the equations
$$V_n(t_i)= \sigma_1^2 C_{H_1}t_i^{H_1+1}+\sigma_2^2C_{H_2}t_i^{H_2+1}, i=1,2.$$
Solving these equations, we obtain the estimators
$$\hat \sigma_{1N}^2=\frac{t_1^{-(H_2+1)}V_n(t_1)-t_2^{-(H_2+1)}V_N(t_2)}{C_{H_1}(t_1^{H_1-H_2}-t_2^{H_1-H_2})}$$
and
$$\hat \sigma_{2N}^2=\frac{t_1^{-(H_1+1)}V_n(t_1)-t_2^{-(H_1+1)}V_N(t_2)}{C_{H_2}(t_1^{H_2-H_1}-t_2^{H_2-H_1})}$$
for $\sigma_1^2$ and for $\sigma_2^2$ respectively. From the almost sure convergence of $V_n(t)$ to $\mu(t)$ as $N \raro \ity,$ it follows that 
$$\hat \sigma_{iN}^2 \raro \sigma_i^2 \;\;a.s \;\;\mbox{as}\;\; N \raro \ity$$
for $i=1,2$ whenever $H_1 \neq H_2.$ Observe that the random vector
$$(\sqrt{N}\hat (\sigma_{1N}^2 -\sigma_1^2),\sqrt{N}(\hat \sigma_{2N}^2 -\sigma2^2))$$ 
is linear function of the random vector 
$$(\sqrt{N} (V_n(t_1)-\mu(t_1)),\sqrt{N} (V_n(t_2)-\mu(t_2))$$
with coefficients depending on $t_1,t_2, H_1$ and $H_2.$ Furthermore the random vector
$$(\sqrt{N} (V_n(t_1)-\mu(t_1)),\sqrt{N} (V_n(t_2)-\mu(t_2))$$
is asympotically bivariate normal with mean zero and suitable covariance matrix. Hence, it follows that the random vector
$$(\sqrt{N}\hat (\sigma_{1N}^2 -\sigma_1^2),\sqrt{N}(\hat \sigma_{2N}^2 -\sigma_2^2))$$
is asymptotically bivariate normal with mean zero and suitable covariance matrix $\Sigma.$
\vsp
Following the method of moments, one can obtain alternate set of estimators by observing that, for any $\delta >o,$ 
$$\frac{1}{N}\sum_{k=1}^N[u(t, k\delta)]^4 \raro  \mu_4 \;\;a.s \;\; \mbox{as}\;\; N \raro \ity$$
where $\mu_4$ is the 4-th central moment of the Gaussian distribution with mean zero and variance $3[\sigma_1^2v_t(H_1)+\sigma_2^2v_t(H_2]^2.$
This follows from the observation that for a Gaussian  distribution with mean zero and variance $\sigma^2$, the 4-th central moment is $3 \sigma^4.$ Let
$$J_N(t) = \frac{1}{N}\sum_{k=1}^N[u(t, k\delta)]^4.$$
\vsp
Solving the equations
$$V_N(t_1)=\sigma_1^2 C_{H_1}t_1^{H_1+1}+\sigma_2^2 C_{H_2}t_1^{H_2+1}$$ 
and
$$\sqrt{J_N(t_2)/3}= \sigma_1^2 C_{H_1}t_2^{H_1+1}+\sigma_2^2 C_{H_2}t_2^{H_2+1},$$
we obtain alternate estimators for $\sigma_1^2$ and $\sigma_2^2$ depending on $H_1,H_2$ and the choice of $t_1$ and $t_2.$ These estimators will also be strongly consistent as $N \raro \ity.$\\
\vsp
\NI{\bf Acknowledgment:} This work was supported under the scheme ``INSA Senior Scientist" at the CR Rao Advanced Institute of Mathematics, Statistics and Computer Science, Hyderabad 500046, India.\\
\vsp
\NI {\bf References:}
\begin{description}
\item Arcones, M. (1994) Limit theorems for nonlinear functionals of a stationary Gaussian sequences of vectors, {\it Ann. Probab.}, {\bf 22}, 2242-2274.
\item Avetisian, D. and Ralchenko, K. (2020) Ergodic properties of the solution to a fractional stochastic heat equation with an application to diffusion parameter estimation, {\it Mod. Stoch. Theory Appl.},{\bf 7}, 339-356.
\item Avetisian, D. and Ralchenko, K. (2023) Parameter estimation in mixed fractional stochastic heat equation,{\it Mod. Stoch. Theory Appl.},{\bf 10}, 175-195.
\item Mishura, Y. and Zili, M. (2018) {\it Stochastic Analysis of Mixed Fractional Gaussian Processes}, ISTE Press and Elsevier, London.
\item Prakasa Rao, B.L.S. (2025) On the long range dependence of time-changed generalized mixed fractional Brownian motion, {\it Calcutta Stat. Assoc. Bull.}, https://doi.org/10.1177/00080683251317814.

\end{description}
\end{document}